\input amstex
\documentstyle{amsppt}

\hcorrection{19mm}

\nologo
\NoBlackBoxes



\topmatter
\title     Incompressible surfaces in link complements
\endtitle
\author    Ying-Qing Wu
\endauthor
\address Department of Mathematics, University of Iowa, Iowa City, IA
52242
\endaddress
\email  wu\@math.uiowa.edu
\endemail

\thanks Supported in part by NSF grant \#DMS 9802558.
\endthanks




\abstract 
We generalize a theorem of Finkelstein and Moriah and show that if a
link $L$ has a $2n$-plat projection satisfying certain conditions,
then its complement contains some closed essential surfaces.  In most
cases these surfaces remain essential after any totally nontrivial
surgery on $L$.
\endabstract

\endtopmatter
 
\document
\define\proof{\demo{Proof}}
\define\endproof{\qed \enddemo}
\redefine\tilde{\widetilde}
\redefine\hat{\widehat}
\define\bdd{\partial}
\define\Int{\text{\rm Int}}
\define\fig#1#2{\leavevmode

                \centerline{\epsfbox{#1}} \medskip  
                \centerline{#2} \medskip }

\input epsf.tex
\TagsOnRight

A link $L$ in $S^3$ has a $2n$-plat projection for some $n$, as shown
in Figure 1, where a box on the $i$-th row and $j$-th column consists
of 2 vertical strings with $a_{ij}$ left-hand half twist; in other
words, it is a rational tangle of slope $1/a_{ij}$.  See for example
[BZ].  Let $n$ be the number of boxes in the even rows, so there are
$n-1$ boxes in the odd rows.  Let $m$ be the number of rows in the
diagram.  It was shown by Finkelstein and Moriah [FM1, FM2] that if
$n\geq 3$, $m \geq 5$, and if $|a_{ij}| \geq 3$ for all $i, j$, then
the link exterior $E(L) = S^3 - \Int N(L)$ contains some essential
planar surfaces, which can be tubed on one side to obtain closed
incompressible surfaces in $E(K)$.  In this note we will prove a
stronger version of this theorem, showing that $E(L)$ contains some
essential surfaces if $n\geq 3$, the boxes on the two ends of the odd
rows have $|a_{ij}| \geq 3$, and $a_{ij} \neq 0$ for the boxes which
are not on the ends of the rows.  We allow $a_{ij} = 0$ for boxes on
the ends of the even rows, and there is no restriction on $m$, the
number of rows in the diagram.  The argument here provides a much
simpler proof to the above theorem of Finkelstein and Moriah.  In
[FM2] that theorem was applied to show that if $L$ is a knot then all
surgeries on $L$ contain essential surfaces.  Corollary 2 below
generalizes this to the case when $L$ has multiple components, with a
mild restriction that each component of $L$ intersects some
``allowable'' spheres.

We first give some definitions.  Let $\alpha = \alpha(a_1, \ldots,
a_m)$ be an arc running monotonically from the top to the bottom of
the $2n$-plat, such that $\alpha$ is disjoint from the boxes, and on
the $i$-th row there are $a_i$ boxes on the left of $\alpha$.  See
Figure 1 for the arc $\alpha(1,1,1,2,2)$.  The arc $\alpha$ is an {\it
allowable path\/} if (i) each row has at least one box on each side of
$\alpha$, and (ii) $\alpha$ intersects $L$ at $m+1$ points, (so
$\alpha$ intersects $L$ once when passing from one row to another).
Note that the leftmost allowable path is $\alpha(1, \ldots, 1)$, which
has on its left one box from each row.

\fig{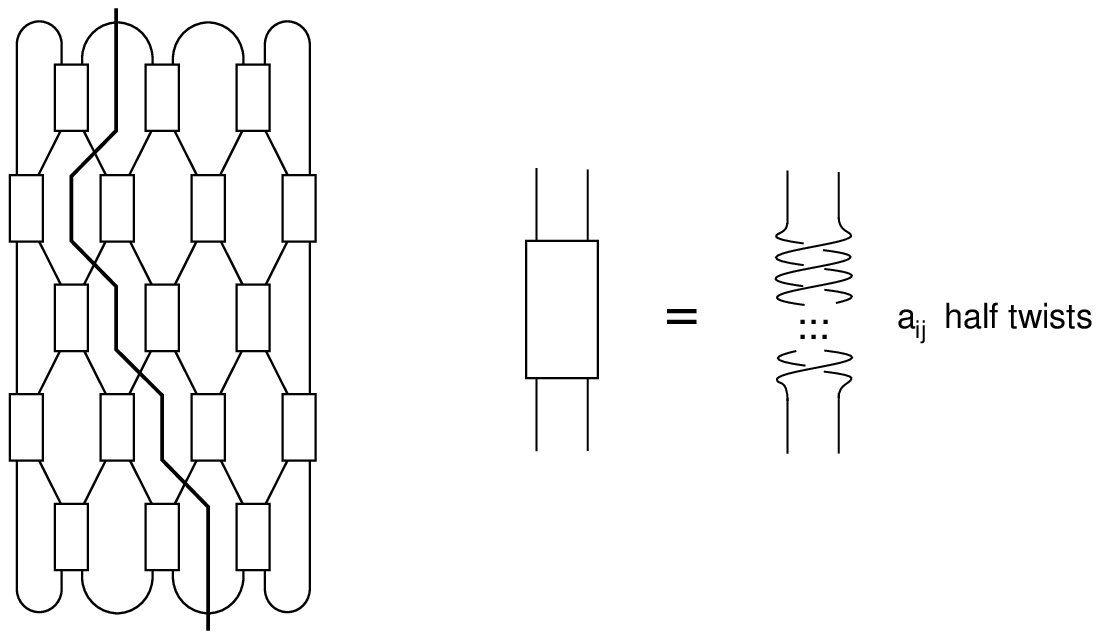}{Figure 1}

Given an allowable path $\alpha = \alpha(a_1, \ldots, a_m)$, we can
connect the two ends of $\alpha$ by an arc $\beta$ disjoint from the
projection of $L$ to form a circle, then cap it off by two disks, one
on each side of the projection plane, to get a sphere $S = S(a_1,
\ldots, a_m)$, called an {\it allowable sphere.}  $S$ cuts $(S^3, L)$
into two tangles $(B, T)$ and $(B', T')$, where $(B, T)$ denotes the
one on the left hand side of $S$.  Let $P = P(a_1, \ldots, a_m)$ be
the planar surface $S \cap E(L)$, which cuts $E(L)$ into two pieces $X
= X(a_1, \ldots, a_m)$ and $X'= X'(a_1, \ldots, a_m)$, with $X = B
\cap E(L)$ the one on the left of $P$.  Let $F = F(a_1, \ldots, a_m)$
be the surface obtained by tubing $P$ on the left hand side; in other
words, $F$ is the component of $\bdd X$ containing $P$, pushed
slightly into the interior of $E(L)$.  Similarly, denote by $F' =
F'(a_1, \ldots, a_m)$ the surface obtained by tubing $P$ on the right
hand side.

Recall that a properly embedded surface $F$ in a 3-manifold $M$ is an
{\it essential surface\/} if it is incompressible,
$\bdd$-incompressible, and is not boundary parallel.  We define a
surface $F$ on the boundary of $M$ to be {\it essential\/} if it is
incompressible, $M \neq F\times I$, and there is no compressing disk
of $\bdd M$ which intersects $F$ at a single essential arc in $F$.
Thus if $F$ is properly embedded in $M$, then it is essential if and
only if after cutting along $F$ the two copies of $F$ are essential in
the resulting manifold.  A 3-manifold $M$ is {\it
$\bdd$-irreducible\/} if $\bdd M$ is incompressible in $M$.  Given a
set $A$ in $M$, denote by $N(A)$ a regular neighborhood of $A$ in $M$.

\proclaim{Theorem 1} Suppose $L$ has a $2n$-plat projection such that
(i) $n\geq 3$; (ii) $a_{ij} \neq 0$ for $j \neq 0, n$; and (iii)
$|a_{ij}| \geq 3$ for $i$ odd and $j=0$ or $n-1$.  Let $S = S(a_1,
\ldots, a_m)$ be an allowable sphere.  Then $E(L)$ is irreducible, and
the surfaces $F= F(a_1, \ldots, a_m)$ and $F'= F'(a_1, \ldots, a_m)$
are essential in $E(L)$.  \endproclaim

Let $L = L_1 \cup \ldots L_k$ be a $k$ component link, let $r = (r_1,
\ldots, r_k)$ be a set of slopes on $\bdd N(L)$, with $r_i$ a slope on
$\bdd N(L_i)$.  Then $L(r)$ denotes the $r$-Dehn surgery on $L$, which
is the manifold obtained by gluing $k$ solid tori $V_1, \ldots,
V_k$ to $E(L)$ so that each $r_i$ is identified with a meridian
disk of $V_i$.  The surgery and the slope $r$ are {\it totally
nontrivial\/} if no $r_i$ is the meridian slope of $L_i$.

\proclaim{Corollary 2}  Let $L$ be as in Theorem 1.  If each component
of $L$ intersects some allowable sphere, then $L(r)$ is a Haken
manifold for all totally nontrivial $r$, and the surfaces $F$ and $F'$
in Theorem 1 remains incompressible in $L(r)$.
\endproclaim

\remark{Remark} (1) It is easy to see that $F=F(a_1, \ldots, a_m)$
being incompressible implies that $P = P(a_1, \ldots, a_m)$ is an
essential planar surface in $E(L)$.  With a similar proof to that of
Theorem 1 one can show that $P$ is essential even if the condition
$|a_{ij}| \geq 3$ in (iii) of Theorem 1 is replaced by $|a_{ij}| \geq
2$.  This generalizes the main theorem of [FM1].

(2) When $n \leq 2$, the link is a 2-bridge link, so by [HT] $E(L)$
contains no closed essential surface.  Hence the assumption 
$n\geq 3$ in Theorem 1 is necessary.

(3) By definition of $2n$-plat projection, the number of rows $m$ is
odd.  If $m = 1$ the link is a composite link, and our assumption
implies that it is nonsplit.  In this case $E(L)$ is irreducible, and
the surfaces in the theorem are swallow-follow tori, which are
essential.  Therefore the theorem is true for $m = 1$.  We may thus
assume that $m \geq 3$ in the proof of Theorem 1.

(4) In Corollary 2, each component of $L$ intersects some allowable
sphere if and only if no component of $L$ is on the left of
$S(1,\ldots, 1)$ or the right of $S(n-2, n-1, \ldots, n-2)$, which is
equivalent to that $a_{i1}$ and $a_{j,n-1}$ are odd for some odd $i,
j$.  

(5) The results remain true if we replace the twist tangles with
rational tangles of slopes $p_{ij}/a_{ij}$ with $a_{ij}$ satisfying
the conditions in the theorem, or certain kinds of more complicated
tangles.  However in this case the link diagram would not be in
$2n$-plat form.  \endremark

\medskip

A $p/q$ {\it rational tangle\/} is a pair $(B, T)$, where $B$ is a
``pillow case'' in ${\bold R}^3$ with corner points $(0, \pm 1, \pm
1)$, and $T$ is obtained by taking 2 arcs of slope $p/q$ on $\bdd B$
connecting the four conner points of the pillow case, then pushing the
interior of the arcs into the interior of $B$.  The $xz$-plane
intersects $\bdd B$ in a circle $C$ of slope $\infty$, called a {\it
vertical circle\/} on $\bdd B$.  Each component of $\bdd B - C$
contains two points of $\bdd T$.  We need the following result about
rational tangles.

\proclaim{Lemma 3} Suppose $(B, T)$ is a $p/q$ rational tangle, and
$C$ a vertical circle on $\bdd B$.  Let $X = B - \Int N(T)$, and let
$P$ be a component of $(\bdd B\cap X) - C$.

(i) if $q \geq 1$ then $P$ is incompressible in $X$; 

(ii) if $q \geq 2$, then $\bdd X - C$ is incompressible in $X$;

(iii) if $q \geq 3$, then any compressing disk of $\bdd X$ intersects
$P$ at least twice.
\endproclaim

\proof
(ii)  Notice that when attaching a 2-handle to $X$ along the curve $C$,
the manifold $X_C$ is the exterior of a 2-bridge link associated to the
rational number $p/q$, which is nontrivial and nonsplit when $q \geq
2$.  In particular, $\bdd X_C$ is incompressible.  If $D$ is a
compressing disk of $\bdd X$ disjoint from $C$, then since $X$ is a
handlebody of genus 2, we can find a nonseparating compressing disk
$D'$ which is still disjoint from $C$.  But then $D'$ would remain a
compressing disk in $X_C$, a contradiction.

(i)  If $q \geq 2$ this follows from (ii) and the fact that $P$ is a
subsurface of $\bdd X - C$ whose complement contains no disk
components.  If $q = 1$, $X$ is a product $P \times I$, and the result
is obvious.

(iii) By (i) $P$ is incompressible, which also implies that $\bdd X -
P$ is incompressible because any simple loop on $\bdd X - P$ is
isotopic to one in $P$.  By [Wu, Lemma 2.1] there is no compressing
disk of $X$ intersecting $P$ at a single essential arc.  \endproof

The following lemma is well-known.  The proof is an easy inner-most
circle outer-most arc argument, and will be omitted.

\proclaim{Lemma 4} Let $F$ be an essential surface in a compact
orientable 3-manifold $M$.  If $M' = M - \Int N(F)$ is irreducible,
and no compressing disk of $\bdd M'$ is disjoint from the two copies
of $F$ on $\bdd M'$, then $M$ is irreducible and $\bdd$-irreducible.
\endproclaim

We now proceed to prove Theorem 1.  In the following, we will assume
that $L$ is a link as in Theorem 1.  By the remark above, we may
assume $m\geq 3$.  

\proclaim{Lemma 5} The manifold $X = X(1, \ldots, 1)$ is irreducible
and $\bdd$-irreducible.  \endproclaim

\proof Consider the tangle $(B, T)$ on the left of $S$.  By an isotopy
of $(B, T)$ we can untwist the boxes in $T$ which lie on the even rows
of the projection of $L$, so the tangle $(B, T)$ is equivalent to the
one shown in Figure 2, where each box corresponds to the first box on
an odd row of the projection of $L$; hence there are $k = (m+1)/2 \geq
2$ boxes, ($k=3$ in Figure 2.)  Let $D_1, \ldots, D_k$ be the disks
represented by the dotted lines in Figure 2, which cuts $(B, T)$ into
$k+1$ subtangles $(B_0, T_0), \ldots, (B_k, T_k)$, where $(B_0, T_0)$
is the one in the middle, which intersects all the $D_i$.  Let $P_i =
D_i \cap X$ be the twice punctured disk in $X$ corresponding to $D_i$.
They cut $X$ into $X_0, \ldots, X_k$, with $X_i = B_i - \Int N(T_i)$
the tangle space of $(B_i, T_i)$.

We want to show that $\cup P_i$ is essential in $X$.  Since each
$(B_i, T_i)$, $i \geq 1$, is a twist tangle with at least 3 twists, by
Lemma 3, the surface $P_i$ is essential in $X_i$.  Now consider $X_0$.
Put $Q =\bdd B_0 - \cup D_i$.  If $D$ is a compressing disk of $Q$ in
$X_0$, then it is a disk in $B_0$ disjoint from $T_0 \cup (\cup D_i)$;
but since $T_0 \cup ( \cup D_i)$ is connected, this would imply that
one side of $D$ is disjoint from all $D_i$, hence $\bdd D$ is a
trivial curve on $Q$, which is a contradiction.  Therefore $Q$ is
incompressible in $X_0$.  Assume there is a disk $D$ in $X_0$ such
that $\bdd D \cap (\cup P_i)$ has only one component.  Since each
string of $T_0$ has ends on different $D_i$, we see that $\bdd D \cap
\bdd N(T_0) = \emptyset$, so $\bdd D \cap (\cup P_i)$ is either a
proper arc in some $D_i$ which separates the two points of $T_0$ on
$D_i$, or it is a circle bounding a disk on $D_i$ containing exactly
one point of $T_0$, or $\bdd D$ can be isotoped into $Q$.  The first
two cases are impossible because then $D$ would be a disk in $B_0$
disjoint from $T_0$ and yet each component of $\bdd B_0 - \bdd D$
contains an odd number of endpoints of $T_0$.  The third case
contradicts the incompressibility of $Q$.  This completes the proof
that $\cup P_i$ is an essential surface in $X$.

Notice that all $X_i$ are handlebodies, and hence irreducible.  Since
$Q$ is incompressible in $X_0$, and by Lemma 3, the surfaces $\bdd X_i
- P_i \subset \bdd X_i - \bdd D_i$ are incompressible in $X_i$ for
$i\geq 1$, it follows from Lemma 4 that $X$ is irreducible and
$\bdd$-irreducible.  \endproof

\fig{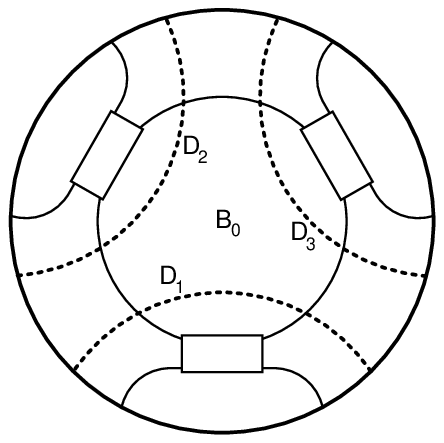}{Figure 2}

\proclaim{Lemma 6} The manifold $X = X(a_1, \ldots, a_m)$ associated
to an allowable sphere $S(a_1, \ldots, a_m)$ is irreducible and
$\bdd$-irreducible.  \endproclaim

\proof There is a sequence of allowable spheres $S_1, \ldots,
S_{k+1}$, such that $S_1 = S(1,\ldots, 1)$, $S_{k+1} = S(a_1, \ldots, a_m)$,
and the non-common part of $S_i, S_{i+1}$ bounds a single box in the
projection of $L$, that is, $S_i \cup S_{i+1} - \Int (S_i \cap
S_{i+1}) = \bdd \tilde B_i$ for some twist tangle $(\tilde B_i, \tilde
T_i)$ with $a \neq 0$ left hand half-twists.  Let $(B_i, T_i)$ be the
tangle on the left of $S_i$, and let $X_i = B_i - \Int N(T_i)$.
Similarly, let $\tilde X_i = \tilde B_i - \Int N(\tilde T_i)$.  Thus
$X = X_{k+1} = X_k \cup _P \tilde X_k$, where $P = X_k \cap
\tilde X_k$ is a twice punctured disk.  By Lemma 5, $X_1$ is
irreducible and $\bdd$-irreducible, and by induction on the length 
of the sequence we may assume that $X_k$ is irreducible and
$\bdd$-irreducible.  Clearly $P$ is incompressible and
$\bdd$-incompressible on the $X_k$ side.  If $|a| \geq 3$ then, by
Lemma 3, $P$ is also incompressible and $\bdd$-incompressible on the
$\tilde X_k$ side, hence $P$ is an essential surface in $X$.
Since $\bdd \tilde X_k - P$ is also incompressible in $\tilde
X_k$, and since $X_k$ and $\tilde X_k$ are irreducible, it
follows that $X = X_k \cup_P \tilde X_k$ is irreducible and
$\bdd$-irreducible.  Also, if $|a| = 1$ then $\tilde X_k$ is a
product $P \times I$, so $X_{k+1} \cong X_k$, and the result follows.

It remains to prove the lemma for the case $|a| = 2$.  In this case
there is a disk $D$ in $\tilde X_k$ which intersects $P$ in a
single arc $\gamma$, cutting $(\tilde X_k, P)$ into a pair
$(A\times I, A \times \bdd I)$, where $A$ is an annulus.  Thus 
$$ X = X_k \cup _P \tilde X_k = (X_k \cup_{\gamma \times I}
(D\times I)) \cup _{A \times \bdd I} ( A\times I) \cong X_k \cup _{A
\times \bdd I} ( A\times I).$$ Since a compressing disk of $\bdd
(A\times I)$ intersects $A\times \bdd I$ at least twice, by the same
argument as above, one can show that $A\times \bdd I$ is essential in
$X$, and $X$ is irreducible and $\bdd$-irreducible.  \endproof

\demo{Proof of Theorem 1} Let $F, F'$ be the surfaces in the theorem,
isotoped slightly to be disjoint from each other.  Then $F\cup F'$
cuts $E(L)$ into three parts: The component on the left of $F$ is
homeomorphic to $X$, the one on the right of $F'$ is homeomorphic to
$X'$, and the one $X''$ between $F$ and $F'$ is the union of $P \times
I$ and $Q\times I$, where $Q$ is the set of tori in $\bdd E(L)$ which
intersect $\bdd P$.  We have shown in Lemma 6 that $X$ is irreducible
and $\bdd$-irreducible, and because of symmetry, so is $X'$.  Now
$X''$ can be cut into $F\times I$ along some (essential) meridional
annuli in $Q \times I$, hence by Lemma 4 it is irreducible and
$\bdd$-irreducible.  Since $F$ and $F'$ have genus at least 2, they
are not boundary parallel.  It follows that $F\cup F'$ is essential in
$X$, and $X$ is irreducible.  \endproof

\demo{Proof of Corollary 2} Let $S_1, \ldots, S_k$ be a set of
disjoint allowable spheres, so that $S_1 = S(1,\ldots, 1)$, $S_k =
S(n-2, n-1, \ldots, n-2)$, and there is only one box of the projection
of $L$ between $S_i$ and $S_{i+1}$.  These spheres are similar to
those in the proof of Lemma 6, except that they are now mutually
disjoint, so the manifold between $S_i$ and $S_{i+1}$ is a product
$S^2 \times I$.  

Let $F_i$ be the essential surfaces corresponding to $S_i$, as defined
before Theorem 1, isotoped slightly so that they are disjoint from
each other.  Also, isotope $F'_k$ to be disjoint from $F_k$.  Then the
set of $k+1$ surfaces $F_1, F_2,\ldots, F_k, F'_k$ cuts $E(L)$ into
$k+2$ components $Y_0, \ldots, Y_{k+1}$, where $Y_0$ is the manifold
$X(1, \ldots, 1)$ on the left of $F_1$, $Y_{k+1}= X'(n-2, n-1, \ldots,
n-2)$ is the manifold on the right of $F_{k+1}$, $Y_k$ is between
$F_k$ and $F_k'$, and for $1\leq i \leq k-1$, $Y_i$ is between $F_i$
and $F_{i+1}$.  Since all the $F_i$ and $F'_k$ are essential, we see
that $Y_i$ are all irreducible and $\bdd$-irreducible.  We need to
show that the manifold $\hat Y_i$ obtained from $Y_i$ by Dehn filling
on its toroidal boundary components (if any), with slopes the
corresponding subset of $r$, is still irreducible and
$\bdd$-irreducible.  The result will then follow by gluing the pieces
together along $F_i$ and $F'_k$.

Our assumption implies that $Y_0$ and $Y_{k+1}$ are disjoint from
$\bdd E(L)$, hence $\hat Y_i = Y_i$ for $i=0, k+1$.  Now $Y_k$ is a
regular neighborhood of $P \cup Q$, where $P = S_k \cap E(L)$, and $Q$
is the set of tori in $\bdd E(L)$ which intersect $P$.  Since $S$ is
separating, each component $Q_j$ of $Q$ intersects $\bdd P$ at least
twice, so there are two nonparallel essential annuli in $Y_k$, each
having a boundary component on $Q_j$ with meridional slope.  Applying
Menasco's theorem [Me] and Scharlemann's theorem [Sch] on each
component of $Q$, we see that after any totally nontrivial Dehn
filling on $Q$ the manifold $\hat Y_k$ is still irreducible and
$\bdd$-irreducible.

Now assume $1\leq i \leq k-1$.  Let $(B'_i, T'_i)$ be the twist tangle
between $S_i$ and $S_{i+1}$.  Notice that if the twist number $a$ of
$T'_i$ is odd then $Y_i$ contains no component of $\bdd E(L)$, so
$\hat Y_i = Y_i$ and we are done.  If $a$ is even, then the tangle
$(B_{i+1}, T_{i+1})$ on the left of $S_{i+1}$ may contain a loop $K$
intersecting the twist tangle $(B'_i, T'_i)$, so $Y_i$ may contain a
single component $Q$ of $\bdd E(L)$.  

\fig{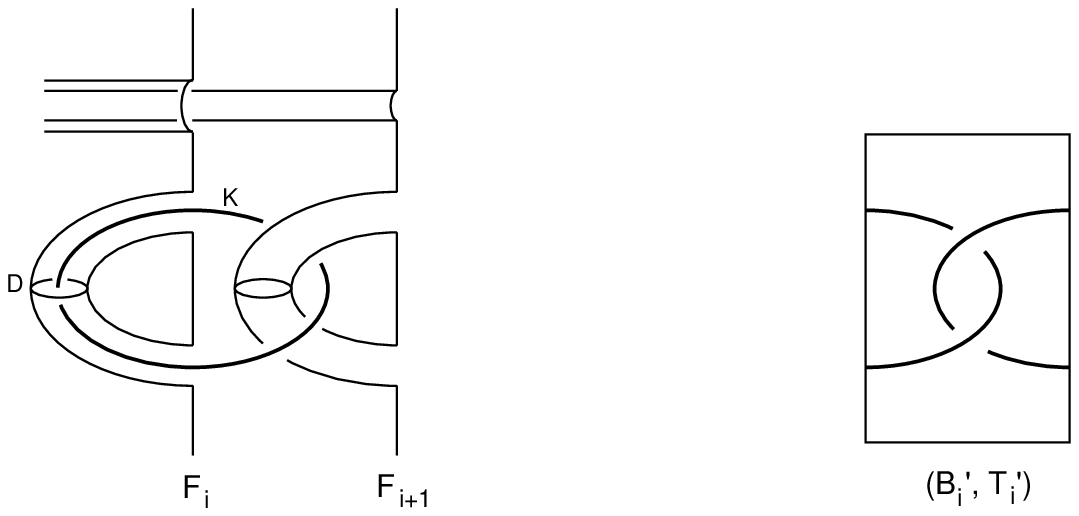}{Figure 3}

Let $Y_i(m)$ be the manifold obtained by the trivial Dehn filling on
$Q$.  Then $F_i$ has a compressing disk $D$ in $Y_i(m)$ intersecting
the core $K$ of the Dehn filling solid torus only once, so $K$ is not
a cable knot in $Y_i(m)$.  See Figure 3.  It follows from [Sch] that
after surgery the manifold $\hat Y_i$ is irreducible.  Also, by [CGLS,
Theorem 2.4.3] $\hat Y_i$ is $\bdd$-irreducible if the surgery slope
$r_j$ on the torus $Q$ intersects the meridian slope $m$ at least
twice.  Now if $r_j$ intersects $m$ only once, then $m$ is a longitude
after the surgery, hence the manifold $\hat Y_i$ is homeomorphic to
the one obtained by cutting $Y_i$ along the annulus $ D \cap Y_i$,
denoted by $\tilde Y_i$.  Now there is an annulus $A$ in $B_{i+1} -
\Int B_i$ ($B_i$ is the ball on the left of $S_i$) separating the
twist tangle $(B'_i, T'_i)$ from the other arcs of $L$, which cuts
$\tilde Y_i$ into $\tilde X\cong B_i' - \Int N(T_i')$ and some
$G\times I$, where $G$ is a subsurface of $F_i$ with one boundary
component.  Clearly $A$ is essential in $G\times I$.  Since the twist
number $a$ is even, our assumption in Theorem 1 implies that $|a| \geq
2$.  Hence by Lemma 3 the surface $\bdd \tilde X - A$ is
incompressible in $\tilde X$, which implies that $A$ is essential in
$\tilde X$.  It follows that $\tilde Y_i$ is irreducible and
$\bdd$-irreducible.  \endproof

\enddocument